\documentclass[11pt]{amsart}
\usepackage{amssymb}

\newtheorem{lemma}{Lemma}[section]

\newtheorem{theorem}[lemma]{Theorem}
\newtheorem{remark}[lemma]{Remark}
\newtheorem{proposition}[lemma]{Proposition}
\newtheorem{corollary}[lemma]{Corollary}
\newtheorem{definition}[lemma]{Definition}

\begin{document}

\title{Euler obstruction and Lipschitz-Killing curvatures}
\author{Nicolas Dutertre }
\address{Aix-Marseille Universit\'e, CNRS, Centrale Marseille, I2M, UMR 7373,
13453 Marseille, France.}
\email{nicolas.dutertre@univ-amu.fr}


\maketitle

\begin{abstract} 
Applying a local Gauss-Bonnet formula for closed subanalytic sets to the complex analytic case, we obtain  characterizations of the Euler obstruction of a complex analytic germ in terms of the Lipschitz-Killing curvatures and the Chern forms of its regular part. 
We also prove analogous results for the global Euler obstruction. As a corollary, we give a positive answer to a question of Fu on the Euler obstruction and the Gauss-Bonnet measure.
 \end{abstract}

\markboth{N. Dutertre}{Euler obstruction and Lipschitz-Killing curvatures}

\section{Introduction}
In \cite{DutertreGeoDedicata} Corollary 5.2, we established the following local Gauss-Bonnet formula for closed subanalytic sets:
$$1=\sum_{k=0}^n    \lim_{\epsilon \rightarrow 0} \frac{\Lambda_k(X,X \cap B_\epsilon)}{b_k \epsilon^k},$$
where $X \subset \mathbb{R}^n$ is a closed subanalytic set containing $0$, $B_\epsilon$ is the closed ball of radius $\epsilon$ centered at $0$, $b_{k}$ is the volume of the unit ball of dimension $k$  and the $\Lambda_k(X,-)$'s are the Lispchitz-Killing measures introduced by Fu in \cite{Fu94}. 

The aim of this paper is to apply this formula to the complex analytic case. In Section 4, we consider a reduced complex analytic germ $(X,0) \subset (\mathbb{C}^n,0)$ of dimension $d>0$ equipped with a finite complex analytic stratification. We write $X= \cup_{i=0}^q V_i$ where $V_0$ is the stratum that contains $0$. 
Let $\phi : X \rightarrow \mathbb{Z}$ be a constructible function. Our main result is the following consequence of the above Gauss-Bonnet formula (see Theorem \ref{MainTheorem}):
  $$\displaylines{
 \quad \phi (0) = \eta(V_0,\phi) + \hfill \cr
\quad \quad \sum_{i=1}^q \left[ \sum_{e=d_0+1}^{d_i} \lim_{\epsilon \rightarrow 0} \frac{1}{b_{2e} \epsilon^{2e}} \frac{1}{s_{2n-2e-1}} \int_{V_i \cap B_\epsilon} K_{2d_i -2e}(x) dx \right] \eta(V_i, \phi) = \hfill\cr
 \hfill  \eta(V_0,\phi) + 
 \sum_{i=1}^q \left[ \sum_{e=d_0+1}^{d_i} \lim_{\epsilon \rightarrow 0} \frac{1}{e ! b_{2e} \epsilon^{2e}} \int_{V_i \cap B_\epsilon}   {\rm ch}_{d_i-e}(V_i) \wedge \kappa(V_i) ^e \right] \eta(V_i, \phi), \qquad \cr 
 }$$
where for each $i \in \{0,\ldots,q\}$, $d_i$ is the dimension of the stratum $V_i$, $K_{2d_i-2e}$ is the $(2d_i-2e)$-th Lipschitz-Killing curvature of $V_i$, ${\rm ch}_{d_i-e}(V_i)$ is the $(d_i-e)$-th Chern form of $V_i$, $\kappa(V_i)$ is the K\"ahler form of $V_i$ and $\eta(V_i,\phi)$ is the normal Morse index of $\phi$ along $V_i$ (see Definition \ref{NMI}). Here $s_{2n-2e-1}$ is the volume of unit sphere of dimension $2n-2e-1$.
Applying this equality to the case where $X$ is equidimensional, we obtain in Corollary \ref{Euler}   characterizations of the Euler obstruction ${\rm Eu}_X(0)$ (see Definition \ref{DefEuler}) in terms of the Lipschitz-Killing curvatures and the Chern forms of the regular part $X_{{\rm reg}}$ of $X$ and in Corollary \ref{IndiceBDK} a new proof of the local index formula of Brylisnki, Dubson and Kashiwara \cite{BDK}. One should mention that, using the relation between the Euler obstruction and the polar multiplicities of L\^e and Teissier \cite{LeTeissierAnnals}, Loeser \cite{Loeser} (p.227) gave a characterization of the Euler obstruction in terms of the Chern-Weil forms of two hermitian bundles, which looks very similar to ours.
We believe that  it should be possible to give a ``complex" proof of our result using Loeser's equality.

 In Section 5, using a Gauss-Bonnet formula for closed semi-algebraic sets (see \cite{DutertreGeoDedicata} Theorem 3.3), we prove global versions of all these local formulas. Here the Euler obstruction is replaced with the global Euler obstruction introduced by Seade, Tib\u{a}r and Verjovsky in \cite{SeadeTibarVerjovsky}.

In \cite{FuJDiffGeo} Section 5, Fu asked a question on the Euler obstruction and the Gauss-Bonnet measure. The Gauss-Bonnet measure is the Lipschitz-Killing measure $\Lambda_0(X,-)$. More precisely, Fu suggested that ${\rm Eu}_X(0)$ should in the limit be equal to the Gauss-Bonnet curvature
of $X \cap B_\epsilon$ within $X_{\rm reg} \cap S_\epsilon$. Then he gave, following in own words, ``two
loosely connected remarks supporting this assertion". In Section 6, using our previous study of the Lipschitz-Killing curvatures in the complex analytic case, we give a proof of Fu's assertion (Corollary \ref{FuQuestion}). 

Throughout the paper, we will use the following notations and conventions (some of them have already appeared in this
introduction):
\begin{itemize}
\item $s_k$ is the volume of unit sphere of dimension $k$ and $b_k$ is the volume of the unit ball of dimension $k$,
\item for $k \in \{0,\ldots,n \}$, $G_{n}^k$ is the Grassmann manifold of $k$-dimension  linear spaces in $\mathbb{R}^{n}$,
$g_n^k$ is its volume,
\item if $H$ is a linear subspace of $\mathbb{R}^{n}$, $S_H$ is the unit sphere in $H$ and $H^\perp$ is the normal space to $H$,
\item for $v \in \mathbb{R}^n$, the function $v^* : \mathbb{R}^n \rightarrow \mathbb{R}$ is defined by $v^*(y)= \langle v, y \rangle$,
\item $B_\epsilon(x)$ is the closed ball of radius $\epsilon$ centered at $x$ and $S_\epsilon(x)$ is the sphere of radius $\epsilon$ centered at $x$, if $x=0$, we simply write $B_\epsilon$ and $S_\epsilon$,
\item if $X \subset \mathbb{R}^{n}$,   $\overline{X}$ is its
topological closure and $\mathring{X}$ is its topological interior.
\end{itemize}

The paper is organized as follows. In Section 2, we introduce the Lipschitz-Killing curvatures  and the polar invariants of Comte and Merle \cite{ComteMerle} and we present results that we proved in  \cite{DutertreGeoDedicata} and \cite{DutertreProcTrotman}. We also give the relation between the Lipschitz-Killing curvatures and the Chern forms. In Section 3, we recall the definitions of the Euler obstruction, the complex link and constructible functions. In Section 4, we apply the local Gauss-Bonnet formula for closed subanalytic sets to the complex case and established our main results. In Section 5, we study the global situation. Section 6 is devoted to the proof of Fu's question.

The reader can refer to \cite{Dubson1}, \cite{Dubson2}, \cite{DutertreGrulha} and \cite{Loeser} for other relations between the Euler obstruction and curvatures.

The author is grateful to Jean-Paul Brasselet and Nivaldo Grulha for introducing him to the subject of Euler obstruction. He is also grateful to J\"org Sch\"urmann for valuable discussions on this topic and for explaining to him the language of constructible functions and to Joe Fu and Andreas Bernig for giving him the relation between the Lipschitz-Killing curvatures and the Chern forms.

\section{The local Gauss-Bonnet formula  for closed subanalytic sets}\label{LocalGB}
In this section, we present the Lipschitz-Killing measures of a subanalytic set and we recall the local Gauss-Bonnet formula for closed subanalytic sets.

In \cite{Fu94}, Fu developped integral geometry for compact subanalytic sets. Using the technology of the normal cycle, he associated with every compact subanalytic set $X$ of $\mathbb{R}^n$ a sequence of curvature measures
$$\Lambda_0(X,-),\ldots,\Lambda_n(X,-),$$
called the Lipschitz-Killing measures. In \cite{BroeckerKuppe} (see also \cite{BernigBroecker}), Broecker and Kuppe gave a geometric characterization of these measures using stratified Morse theory. Let us describe their work.

Let $X \subset \mathbb{R}^n$ be a compact subanalytic set equipped with a finite subanalytic stratification $\mathcal{V}=\{ V_a \}_{a \in A}$. Let us assume that $X$ has dimension $d$ and let us fix a stratum $V$ of dimension $e$ with $e<n$.
%For $x \in X$, we denote by $S_{T_x V^\perp}$ the unit sphere of $T_x V^\perp$, the normal space to $V$ at $x$. 
For $k \in \{0,\ldots,e\}$ and for $x \in X$, let $\lambda_k^V(x) : V \rightarrow \mathbb{R}$ be defined by
$$ \lambda_k^V(x) = \frac{1}{s_{n-k-1}} \int_{S_{T_x V^\perp} }\alpha(x,v) \sigma_{e-k} (II_{x,v}) dv,$$
where $II_{x,v}$ is the second fundamental form on $V$ in the direction of $v$ and $\sigma_{e-k} (II_{x,v})$ is the $(e-k)$-th elementary symmetric function of its eigenvalues. We recall that $II_{x,v}$ is defined on $T_x V$ as follows:
$$II_{x,v}(W_1,W_2) = -  \langle \nabla_{W_1} W, W_2 \rangle,$$
for $W_1$ and $W_2$ in $T_x V$, where $\nabla$ is the covariant differentiation in $\mathbb{R}^n$ and $W$ is a local extension of $v$ normal to $V$. The index $\alpha(x,v)$ is defined as follows:
$$\alpha(x,v)= 1-\chi \Big( X \cap N_x \cap B_\epsilon(x) \cap \{ v^*= v^*(x)-\delta \} \Big),$$
where $0 < \delta \ll \epsilon \ll 1$ and $N_x$ is a normal (subanalytic) slice to $V$ at $x$ in $\mathbb{R}^n$ such that $N_x \cap V = \{x\}$. % $B_\epsilon (x)$ is the closed ball centered at $x$ of radius $\epsilon$ in $\mathbb{R}^n$ and the function $v^* : \mathbb{R}^n \rightarrow \mathbb{R}$ is given by $v^*(y)=\langle v, y \rangle$. 
Since we work in the subanalytic setting, this index is well-defined thanks to Hardt's theorem \cite{Hardt}. Furthermore when $v^*_{\vert X}$ has a stratified Morse critical point at $x$, it coincides with the normal Morse index at $x$ of a function $f : \mathbb{R}^n \rightarrow \mathbb{R}$ such that $f_{\vert X}$ has a stratified Morse critical point at $x$ and $\nabla f (x) = v$ (see \cite{GoreskyMacPherson}, I.1.8 and \cite{BroeckerKuppe}, Lemma 3.5). For $k \in \{e+1,\ldots,n\}$, we set $\lambda_k^V(x)=0$. 

If $V$ has dimension $n$ then for all $x \in V$, we put $\lambda_0^V(x)=\cdots=\lambda_{n-1}^V(x)=0$ and $\lambda_n^V(x)=1$. If $V$ has dimension $0$ then we set
$$\lambda_0^V(x)= \frac{1}{s_{n-1}} \int_{S^{n-1}}  \alpha(x,v) dv,$$
and $\lambda_k^V(x)=0$ if $k>0$. 
\begin{definition}
For every Borel set $U \subset X$ and for every $k \in \{0,\ldots,n\}$, we define $\Lambda_k(X,U)$ by
$$\Lambda_k(X,U)= \sum_{a \in A} \int_{V_a \cap U} \lambda_k^{V_a} (x) dx.$$
\end{definition}
These measures $\Lambda_k(X,-)$ are called the Lipschitz-Killing measures of $X$. Note that for any Borel set $U$ of $X$, we have
$$\Lambda_{d+1}(X,U)= \cdots=\Lambda_n(X,U)=0,$$ 
and $\Lambda_d(X,U)= \mathcal{L}_d(U)$, where $\mathcal{L}_d$ is the $d$-th dimensional Lebesgue measure in $\mathbb{R}^n$. 
If $X$ is smooth then for $k \in \{0,\ldots,d \}$, $\Lambda_k(X,U)$ is equal to 
$$\frac{1}{s_{n-k-1} } \int_U K_{d-k}(x) dx,$$
where $K_{d-k}$ is the $(d-k)$-th Lipschitz-Killing-Weyl curvature. 

Furthermore when $X$ is a smooth complex analytic subset of dimension $d$ in $\mathbb{C}^n$, the Lipschitz-Killing measures are related to the Chern forms of $X$. 
\begin{lemma}\label{LipschitzChern}
Let $X \subset \mathbb{C}^n$ be a smooth  complex analytic subset of dimension $d$. For any Borel set $U$ of $X$ and for $k \in \{0,\ldots, d\}$, we have
$$\Lambda_{2k}(X,U)= \frac{1}{s_{2n-2k-1} } \int_U K_{2d-2k}(x) dx =\frac{1}{k !} \int_U {\rm ch}_{d-k}(X) \wedge \kappa(X)^k,$$
where ${\rm ch}_i (X)$ is the $i$-th Chern form and $\kappa(X)$ is the K\"ahler form of $X$.
\end{lemma}
\proof Apply (2.7.1) and Lemma 2.33 in \cite{BernigFuSolanes} (or Theorem 4.8 and Lemma 7.6 in \cite{Gray}). \endproof

In \cite{DutertreGeoDedicata} Section 5, we studied the asymptotic behavior of the Lipschitz-Killing measures in the neighborhood of a point of $X$. We considered a  closed subanalytic set $X$ such that $0 \in X$ and we showed the following theorem (\cite{DutertreGeoDedicata}, Theorem 5.1):
\begin{theorem}\label{CurvAndLink}
Let $X \subset \mathbb{R}^n$ be a  closed subanalytic set such that $0 \in X$. 
We have
$$\lim_{\epsilon \rightarrow 0} \Lambda_0(X,X\cap B_\epsilon)=1-\frac{1}{2} \chi (\hbox{\em Lk}(X))-\frac{1}{2g_n^{n-1}} \int_{G_n^{n-1}} \chi (\hbox{\em Lk}(X \cap H)) dH.$$
Furthermore for $k \in \{1,\ldots,n-2 \}$, we have
$$\displaylines{
\qquad \lim_{\epsilon \rightarrow 0} \frac{\Lambda_k(X,X\cap B_\epsilon)}{b_k \epsilon^k} =- \frac{1}{2g_n^{n-k-1}} \int_{G_n^{n-k-1}} \chi(\hbox{\em Lk}(X \cap H)) dH \hfill \cr
\hfill + \frac{1}{2g_n^{n-k+1}} \int_{G_n^{n-k+1}} \chi(\hbox{\em Lk}(X \cap L)) dL, \qquad  \cr
}$$
and:
$$\lim_{\epsilon \rightarrow 0} \frac{\Lambda_{n-1}(X,X\cap B_\epsilon)}{b_{n-1} \epsilon^{n-1}} = \frac{1}{2g_n^2} \int_{G_n^2} \chi(\hbox{\em Lk}(X \cap H)) dH,$$
$$\lim_{\epsilon \rightarrow 0} \frac{\Lambda_{n}(X,X\cap B_\epsilon)}{b_{n} \epsilon^{n}} = \frac{1}{2g_n^1} \int_{G_n^1} \chi(\hbox{\em Lk}(X \cap H)) dH.$$
\end{theorem}
$\hfill \Box$

As a corollary, we obtained (\cite{DutertreGeoDedicata}, Corollary 5.2):
\begin{corollary}\label{CorCurvAndLink}
We have
$$1=\sum_{k=0}^n    \lim_{\epsilon \rightarrow 0} \frac{\Lambda_k(X,X \cap B_\epsilon)}{b_k \epsilon^k}.$$
\end{corollary}
$\hfill \Box$

In \cite{DutertreProcTrotman}, we continued our study of the Lipschitz-Killing measures and related the above limits $\lim_{\epsilon \rightarrow 0} \frac{\Lambda_k(X,X \cap B_\epsilon)}{b_k \epsilon^k}$ to the polar invariants introduced by Comte et Merle in \cite{ComteMerle}. These invariants can be defined as follows. Let $H \in G_n^{n-k}$, $k \in \{1,\ldots, n \}$, and let $v$ be an element in $S_{H^\perp}$. For $\delta >0$, we denote by $H_{v,\delta}$ the $(n-k)$-dimensional affine space $H+\delta v$ and we set 
$$\alpha_0(H,v) = \lim_{\epsilon \rightarrow 0} \lim_{\delta \rightarrow 0} \chi(H_{\delta,v} \cap X \cap B_\epsilon),$$
$$\alpha_0(H)=\frac{1}{s_{k-1} }\int_{S_H^\perp} \alpha_0 (H,v) dv,$$
and then 
$$\sigma_k (X,0)=\frac{1}{g_n^{n-k}} \int_{G_n^{n-k}} \alpha_0 (H) dH.$$
Moreover, we put $\sigma_0 (X,0)=1$. In \cite{DutertreProcTrotman}, we proved the following theorem.
\begin{theorem}\label{CurvAndPolar}
For $k \in \{0,\ldots,n-1\}$, we have
$$\lim_{\epsilon \rightarrow 0} \frac{\Lambda_k(X,X \cap B_\epsilon)}{b_k \epsilon^k} = \sigma_k(X,0) -\sigma_{k+1}(X,0).$$
Furthermore, we have
$$\lim_{\epsilon \rightarrow 0} \frac{\Lambda_n(X,X \cap B_\epsilon)}{b_n \epsilon^n} = \sigma_n(X,0).$$
\end{theorem}
$\hfill \Box$

Let $d_0$ be the dimension of the stratum that contains $0$. It is explained in \cite{ComteMerle} that $\sigma_k(X,0)=1$ if $0 \le k \le d_0$. Therefore, we find that
$$\lim_{\epsilon \rightarrow 0} \frac{\Lambda_{d_0}(X,X \cap B_\epsilon)}{b_{d_0} \epsilon^{d_0}} = 1-\sigma_{d_0+1}(X,0),$$
and for $k < d_0$,
$$\lim_{\epsilon \rightarrow 0} \frac{\Lambda_k(X,X \cap B_\epsilon)}{b_k \epsilon^k} =0.$$
We can refine the above corollary.
\begin{corollary}\label{CorCurvAndLink2}
Let $X$ be a closed subanalytic set of dimension $d$ such that $0 \in X$. Let $d_0$ be the dimension of the stratum that contains $0$. We have
$$1=\sum_{k=d_0}^d    \lim_{\epsilon \rightarrow 0} \frac{\Lambda_k(X,X \cap B_\epsilon)}{b_k \epsilon^k}.$$
\end{corollary}
$\hfill \Box$

\section{Euler obstruction, complex link and constructible functions}

\subsection{The Euler obstruction}

The Euler obstruction was defined by Mac\-Pherson \cite{MacPhersonAnnals74} as a tool to prove the conjecture about existence and unicity of the Chern classes in the singular case. For an overview about the Euler obstruction see \cite{Brasselet,BrasseletGrulha}. Let us now introduce some objects in order to define the Euler obstruction.

Let $(X,0) \subset (\mathbb{C}^n,0)$ be an equidimensional reduced complex analytic germ of dimension $d$ in an open set $U \subset \mathbb{C}^n$. We consider a complex analytic Whitney stratification $\mathcal{V} = \{V_a\}_{a \in A}$ of  $X$. We choose a small representative of $(X,0)$ such that $0$ belongs to the closure of all the strata. We will denote it by $X$ and we will write $X= \cup_{i=0}^q V_i$. We will assume that the strata $V_0,\ldots,V_{q}$ are connected and that the analytic sets $\overline{V_0},\ldots,\overline{V_{q}}$ are reduced. 

Let $G(d,n)$ denote the Grassmanian of complex $d$-planes
in ${\mathbb C}^n$. On the regular part $X_{\rm reg}$ of $X$ the Gauss map
$\phi : X_{ \rm reg} \to U\times G(d,n)$ is well defined by $\phi(x) =
(x,T_x(X_{ \rm reg}))$.

\begin{definition}
The {\it Nash transformation} (or {\it Nash blow up}) $\widetilde X$ of $X$
is the closure of the image {\rm Im}$(\phi)$ in $ U\times G(d,n)$. It is a
(usually singular) complex analytic space endowed with an analytic
projection map $\nu : \widetilde   X \to X$
which is a biholomorphism away from $\nu^{-1}({{\rm Sing}}(X))\,.$
\end{definition}

The fiber of  the tautological bundle  ${\mathcal T}$ over $G(d,n)$, at the point $P \in G(d,n)$,
is the set of the vectors $v$ in the $d$-plane $P$. We still denote by ${\mathcal T}$ the
corresponding trivial extension bundle over $ U \times G(d,n)$. Let $\widetilde
T$ be the restriction of  ${\mathcal T}$ to $\widetilde   X$, with
projection map $\pi$.  The bundle $\widetilde   T$ on $\widetilde
X$ is called {\it the  Nash bundle} of $X$.

An element of $\widetilde   T$ is written $(x,P,v)$ where $x\in U$,
$P$ is a $d$-plane in ${\mathbb C}^n$ based at $x$ and $v$ is a vector in
$P$. We have the following diagram:
$$
\begin{matrix}
\widetilde   T & \hookrightarrow & {\mathcal T} \cr
{\pi} \downarrow & & \downarrow \cr
\widetilde   X & \hookrightarrow & U \times G(d,n) \cr
{\nu}\downarrow & & \downarrow \cr
X & \hookrightarrow & U. \cr
\end{matrix}
$$

Let us recall the original definition of the Euler obstruction, due to Mac\-Pherson \cite{MacPhersonAnnals74}.
Let $z=(z_1, \ldots,z_n)$ be local coordinates in ${\mathbb C}^n$ around $\{ 0 \}$, such that
$z_i(0)=0$. %We denote by $ B_\epsilon   $ and $S_\epsilon   $ the ball and the sphere
%centered at $\{ 0 \}$ and of radius $\epsilon   $ in ${\mathbb C}^N$. 
Let us consider the norm
$\Vert z \Vert = \sqrt{z_1 \overline z_1 + \cdots + z_n \overline z_n}$. Then the differential form
$\omega = d \Vert z \Vert ^2$ defines a section of the real vector bundle $T^{*}{\mathbb C}^n$,
cotangent bundle on ${\mathbb C}^n$. Its pull back restricted to $ \widetilde   X$ becomes a section
denoted by $ \widetilde\omega$ of the dual bundle ${\widetilde   T}^*$. For $\epsilon   $ small enough,
the section  $ \widetilde\omega$ is nonzero over $\nu^{-1}(z)$ for $0 < \Vert z \Vert \le \epsilon   $.
The obstruction to extend $ \widetilde\omega$ as a non-zero section of ${\widetilde   T}^*$ from
$\nu^{-1}(S_\epsilon   )$ to $\nu^{-1}(B_\epsilon   )$, denoted by
$Obs({\widetilde   T}^*, \widetilde\omega)$, lies in $H^{2d}(\nu^{-1}( B_\epsilon   ),
\nu^{-1}( S_\epsilon   );{\mathbb Z})$. Let us denote by ${\mathcal O}_{\nu^{-1}(B_\epsilon   ),
\nu^{-1}( S_\epsilon   )}$ the orientation class in $H_{2d}(\nu^{-1}(B_\epsilon   ),
\nu^{-1}(S_\epsilon   );{\mathbb Z})$.

\begin{definition}\label{DefEuler}
The local Euler obstruction of $X$ at $0$ is the evaluation of $Obs({\widetilde   T}^*, \widetilde\omega)$ on
${\mathcal O}_{\nu^{-1}(\mathbb B_\epsilon   ), \nu^{-1}(\mathbb S_\epsilon   )}$, i.e.:
$${\rm Eu}_{X}(0) = \langle  Obs({\widetilde   T}^*, \widetilde\omega), {\mathcal O}_{\nu^{-1}(B_\epsilon   ),
\nu^{-1}(S_\epsilon   )} \rangle.$$
\end{definition}

An equivalent definition of the Euler obstruction was given by Brasselet and Schwartz in the context of vector fields \cite {BrasseletSchwartz}.

\subsection{The complex link and the normal Morse datum}
The complex link is an important object in the study of the topology of complex analytic sets. It is analogous to the Milnor fibre and was studied first in \cite{Le75}. It plays a crucial role in complex stratified Morse theory (see \cite{GoreskyMacPherson}) and appears in general bouquet theorems for the Milnor fibre of a function with isolated singularity (see \cite{Le86, Le92,SiersmaBouquet, TibarBouquet}).  

We recall first the definitions of a conormal covector, a degenerate tangent plane, a degenerate covector and an exceptional point.
Let $X \subset \mathbb{C}^n$ be  a reduced complex analytic set of dimension $d$. We assume that $X$ is included in an open set $U$ and that is is equipped with a Whitney stratification $\mathcal{V}=\{ V_a \}_{a \in A}$, which strata are connected. 
\begin{definition}
Let $x$ be a point in $X$ and let $V_b$ be the stratum that contains it. A cotangent vector $\eta \in T^*_x U$ is conormal for $X$ at $x$ if $\eta(T_x V_b)=0$.
\end{definition}

\begin{definition}
Let $x$ be a point in $X$ and let $V_b$ be the stratum that contains it. A degenerate tangent plane of the stratification $\mathcal{V}$ at $x$ is an element $T$ (of an appropriate Grassmanian) such that  $T=\lim_{x_i \rightarrow x} T_{x_{i}} V_a$, where $V_a$ is a stratum that contains  $V_b$ in its frontier and where the $x_{i}$'s belong to  $V_a$.
\end{definition}

\begin{definition}
A degenerate covector of $\mathcal{V}$ at a point $x \in X$ is a covector which vanishes on a degenerate tangent plane of $\mathcal{V}$ at $x$, {\it i.e.} an element $\eta$ of $T^{*}_{p}U$ such that there exists a degenerate tangent plane $T$ of $\mathcal{V}$ at $x$ with ${\eta}(T)=0$. 
\end{definition}

\begin{definition}
A point $x$ in $X$ is exceptional if the degenerate conormal vectors at $x$ form a codimension $0$ subvariety of the conormal space at $x$.
\end{definition}

Teissier \cite{TeissierLaRabida} (Prop. 1.2.1, p. 461) proved that a Whitney stratified complex analytic set does not admit exceptional points. 

We can now define the complex link. 
Let $V$ be a stratum of the stratification $\mathcal{V}$ of $X$ and let $x$ be a point in $V$. Let $g : (\mathbb{C}^n,x) \rightarrow (\mathbb{C},0)$ be an analytic complex function-germ such that the differential form $Dg(x)$ is not a degenerate covector of $\mathcal{V}$ at $x$. Let $N^{\mathbb{C}}_x$ be a normal slice to $V$ at $x$, i.e. $N^{\mathbb{C}}_x$ is a closed complex submanifold of $\mathbb{C}^n$ which is transversal to $V$ at $x$ and $N^{\mathbb{C}}_x \cap V =\{x\}$.
\begin{definition}
The complex link $\mathcal{L}_V$ of $V$ is defined by
$$\mathcal{L}_V = X\cap N_x^{\mathbb{C}}  \cap B_{\epsilon}(x)\cap \{g=\delta\} ,$$
where $ 0< \vert \delta \vert \ll \epsilon \ll 1$.  %Here $B_{\epsilon}(x)$ is the closed ball of radius $\epsilon$ centered at $x$ .

The normal Morse datum ${\rm NMD}(V)$ of $V$ is the pair of spaces
$${\rm NMD}(V) =(X\cap N_x^{\mathbb{C}} \cap B_{\epsilon}(x), X\cap N_x^{\mathbb{C}} \cap B_{\epsilon}(x)\cap \{g=\delta\}).$$
\end{definition}
The fact that these two notions are well-defined, i.e. independent of all the choices made to define them, is explained in \cite{GoreskyMacPherson}.  

The set $X$ can be viewed as a real analytic set in $\mathbb{R}^{2n}$ and we can compare $\chi(\mathcal{L}_V)$ to the indices $\alpha(x,v)$, $x\in V$ and $v \in S_{T_x V^\perp}$, introduced in Section 2.  Since ${\rm Hom}_{\mathbb{R}} (T_x V, \mathbb{R})$ is canonically isomorphic to ${\rm Hom}_{\mathbb{C}} (T_x V, \mathbb{C})$, the point $x$ is not exceptional in the real sense. Therefore, for almost all $v$ in $T_x V ^\perp$, the form $v^*$ is a non-degenerate conormal covector at $x$ in the real sense. Hence, for almost all $v$ in $T_x V ^\perp$, the set 
$$ X \cap N_x \cap B_\epsilon(x) \cap \{ v^* = v^*(x)-\delta \},$$
introduced in Section 2 is the lower half-link of the function $v^*_{\vert X}$ and by \cite{GoreskyMacPherson} (Theorem I 3.9.3 and Corollary 1, Section 2.5, Part II), its Euler characteristic is equal to $\chi(\mathcal{L}_V)$. We can conclude that for all $x$ in $V$ and almost all $v$ in $T_x V^\perp$, we have 
$$\alpha(x,v)=1-\chi(\mathcal{L}_V).$$

\subsection{Constructible functions}
Let $(X,0) \subset (\mathbb{C}^n,0)$ be a reduced complex analytic germ of dimension $d$ in a open set $U \subset \mathbb{C}^n$. We consider a complex analytic Whitney stratification $\mathcal{V}=\{V_a\}_{a \in A}$ of  $X$. We choose  a small representative of $(X,0)$ such that $0$ belongs to the closure of all the strata. We still denote it by $X$ and we write $X=\cup_{i=0}^q V_i$. We assume that the stratum $V_0,\ldots,V_q$ are connected and that the analytic sets $\overline{V_0},\ldots,,\overline{V_q}$ are reduced. 
\begin{definition}
A constructible function with respect to the stratification $\mathcal{V}$ of $X$ is a function $\beta: X \to \mathbb{Z}$ which is constant on each stratum $V_i$. 
\end{definition}
This means that there exist integers $n_0,\ldots,n_q$ such that  we can write
$$\beta=\sum_{i=0}^{q}n_{i} \cdot {\bf 1}_{V_{i}},$$
where ${\bf 1}_{V_{i}}$ is the characteristic function on $V_{i}$.

\begin{remark}
When $X$ is equidimensional, there are two distinguished bases for the free abelian group of such constructible functions: the characteristic functions ${\bf 1}_{\overline{V}}$ and the Euler obstruction ${\rm Eu}_{\overline{V}}$ of the closure $\overline{V}$ of all strata $V$.
\end{remark}

\begin{definition}
The Euler characteristic $\chi(X,\beta)$ of a constructible function $\beta : X \rightarrow \mathbb{Z}$ given by $\beta= \sum_{i=0}^q n_i  {\bf 1}_{V_i}$ is defined by
$$\chi(X,\beta) = \sum_{i=0}^q n_i \  \chi(V_i).$$
\end{definition}

\begin{definition}\label{NMI}
Let $\beta : X \to Z$ be a constructible function with respect to the stratification $\mathcal{V}$. Its normal Morse index $\eta(V,\beta)$ along $V$ is defined by
$$\eta(V,\beta)=\chi({\rm NMD}(V),\beta )=\chi(X \cap N \cap B_{\epsilon}(x),\beta)-\chi(\mathcal{L}_{V},\beta).$$
\end{definition}

If $Z\subset X$ is a closed union of strata, then $\eta(V,{\bf 1}_{Z})=1-\chi(\mathcal{L}_{V}\cap Z)$. The key role of the Euler
obstruction comes from the following identities (see \cite{SchuermannTibar} p.34 or \cite{Schuermann} p.292 and p.323-324):
$$\eta(V',{\rm Eu}_{\overline{V}})= 1\hbox{ if }V'=V,$$ and
$$\eta(V', {\rm Eu}_{\overline{V}})=0 \hbox{ if } V' \neq V.$$

\section{Euler obstruction and curvatures}
We apply the result of Section 2 to the case of a complex analytic germ.

Let $(X,0) \subset (\mathbb{C}^n,0)$ be a reduced complex analytic germ of dimension $d>0$ in a open set $U \subset \mathbb{C}^n$ and let $\mathcal{V}=\{V_a\}_{a \in A}$ be a complex analytic Whitney stratification of $X$. We choose  a small representative of $(X,0)$ such that $0$ belongs to the closure of all the strata. We still denote it by $X$ and we write $X=\cup_{i=0}^q V_i$ where $V_0$ is the stratum containing $0$. We assume that the stratum $V_0,\ldots,V_q$ are connected and that the analytic sets $\overline{V_0},\ldots,\overline{V_q}$ are reduced. We observe that $V_0=\overline{V_0}$ as analytic germs. We set $d_i= {\rm dim} V_i$ for $i \in \{0,\ldots,q\}$. We call a {\bf top stratum} a stratum not contained in the frontier of any other stratum. 

Viewing $\mathbb{C}^n$ as $\mathbb{R}^{2n}$ and applying the results of Section 2, we already know that 
$$ \lim_{\epsilon \rightarrow 0} \frac{\Lambda_k(X,X \cap B_\epsilon)}{b_k \epsilon^k}=0,$$
for $k=0,\ldots, 2d_0-1$ and $k=2d+1,\ldots,2n$. If $k$ is an even integer between $2d_0$ and $2d$ and $H$ is a generic element of the Grassmaniann $G_{2n}^{2n-k}$, then ${\rm Lk}(X \cap H)$ is a Whitney stratified set with only odd-dimensional strata so, by Sullivan's theorem \cite{Sullivan}, $\chi({\rm Lk}(X \cap H))=0$. Hence by Theorem \ref{CurvAndLink}, we see that 
$$ \lim_{\epsilon \rightarrow 0} \frac{\Lambda_k(X,X \cap B_\epsilon)}{b_k \epsilon^k}=0,$$
if $k$ is an odd integer between $2d_0$ and $2d$, and so
$$1=\sum_{e=d_0}^d    \lim_{\epsilon \rightarrow 0} \frac{\Lambda_{2e}(X,X \cap B_\epsilon)}{b_{2e} \epsilon^{2e}}.$$
Let us focus now on the polar invariants $\sigma_k(X,0)$. By \cite{ComteMerle}, we know that
$$\sigma_0(X,0)=\cdots=\sigma_{2d_0}(X,0)=1 \hbox{ and } \sigma_{2d+1}(X,0)=\cdots=\sigma_{2n}(X,0)=0.$$
By the above equality and Theorem \ref{CurvAndPolar}, we get that for $e \in \{d_0+1,\ldots,d\}$, $\sigma_{2e-1}=\sigma_{2e}$. Therefore, Theorem \ref{CurvAndPolar} can be rewritten as follows.
\begin{theorem}\label{CurvAndPolarC}
We have 
$$ \lim_{\epsilon \rightarrow 0} \frac{\Lambda_{2d_0} (X,X \cap B_\epsilon)}{b_{2d_0} \epsilon^{2d_0}} = 1-\sigma_{2d_0 +1}(X,0)= 1-\sigma_{2(d_0+1)} (X,0),$$
and for $e \in \{d_0+1,\ldots,d-1\}$
$$\lim_{\epsilon \rightarrow 0} \frac{\Lambda_{2e} (X,X \cap B_\epsilon)}{b_{2e} \epsilon^{2e}} = \sigma_{2e}(X,0)-\sigma_{2e +1}(X,0)= \sigma_{2e}(X,0)-\sigma_{2(e+1)} (X,0),$$
and
$$\lim_{\epsilon \rightarrow 0} \frac{\Lambda_{2d} (X,X \cap B_\epsilon)}{b_{2d} \epsilon^{2d}} = \sigma_{2d}(X,0).$$
\end{theorem}
$\hfill \Box$

Let us study the Lipschitz-Killing curvatures $\Lambda_{2e}(X,X \cap B_\epsilon)$ more carefully. We have
$$\Lambda_{2e}(X,X \cap B_\epsilon) = \sum_{i=0}^q \int_{V_i \cap B_\epsilon} \lambda_{2e}^{V_i} (x) dx,$$
where 
$$\lambda_{2e}^{V_i}(x)= \frac{1}{s_{2n-2e-1}} \int_{S_{T_x V_i ^\perp}} \sigma_{2(d_i-e)} (II_{x,v}) \alpha(x,v) dv,$$ 
if $e \le d_i$ and $\lambda_{2e}^{V_i}(x)=0$ if $e > d_i$. 
But we know that for $x \in V_i$ and for almost all $v $ in $S_{T_x V_i^\perp}$,
$$\alpha(x,v)= 1-\chi(\mathcal{L}_{V_i})=\eta (V_i, {\bf 1}_X).$$
Therefore, for $e \le d_i$, we have
$$\displaylines{
\qquad  \lambda_{2e}^{V_i}(x) = \eta(V_i,{\bf 1}_X) \frac{1}{s_{2n-2e-1}} \int_{S_{T_x V_i ^\perp}} \sigma_{2(d_i-e)} (II_{x,v} ) dv = \hfill \cr
\hfill  \eta(V_i,{\bf 1}_X) \frac{1}{s_{2n-2e-1}}  K_{2(d_i-e)} (x). \qquad \cr
}$$
Putting $K_{2(d_i-e)}(x)=0$ if $e > d_i$, we can write
$$\Lambda_{2e}(X,X\cap B_\epsilon) = \sum_{i=0}^q \frac{\eta(V_i,{\bf 1}_X)}{s_{2n-2e-1}} \int_{V_i \cap B_\epsilon} K_{2(d_i-e)}(x)dx.$$
\begin{proposition}
For $i \in \{0,\ldots,q\}$, for $e \in \{0,\ldots,n\}$, 
$$\lim_{\epsilon \rightarrow 0} \frac{1}{b_{2e} \epsilon^{2e}} \frac{1}{s_{2n-2e-1}}  \int_{V_i \cap B_\epsilon} K_{2(d_i-e)}(x) dx,$$
is finite. Furthermore, for $e \in \{0,\ldots,d_0-1\} \cup \{d_i+1,\ldots,n\}$, this limit vanishes.
\end{proposition}
\proof It is clear that $$\lim_{\epsilon \rightarrow 0} \frac{1}{b_{2e} \epsilon^{2e}} \frac{1}{s_{2n-2e-1}}  \int_{V_i \cap B_\epsilon} K_{2(d_i-e)} (x) dx=0$$ if $e >d_i$ by definition of $K_{2(d_i-e)}$. 

We prove the finiteness of the other limits by induction on the depth of the stratum. Let us start with the stratum $V_0$. In this case, the finitess of the limits is guaranteed by Theorem \ref{CurvAndLink} and, by the remark before Corollary \ref{CorCurvAndLink2}, we have
$$\lim_{\epsilon \rightarrow 0} \frac{1}{b_{2e} \epsilon^{2e}} \frac{1}{s_{2n-2e-1}} \int_{V_0 \cap B_\epsilon} K_{2(d_0-e)} (x) dx =0,$$
for $e \in \{0,\ldots,d_0-1\}$. Furthermore, Corollary \ref{CorCurvAndLink2} implies that
$$\lim_{\epsilon \rightarrow 0} \frac{1}{b_{2 d_0} \epsilon^{2d_0}} \frac{1}{s_{2n-2d_0-1}} \int_{V_0 \cap B_\epsilon} K_0(x)dx =1.$$
Let us fix now a stratum $V_i$. By induction, we can assume that the result is true for any stratum $W$ included in $\overline{V_i} \setminus V_i$. For $e \in \{0,\ldots,d_0-1\}$, we have
$$\displaylines{ 
\qquad \Lambda_{2e}(\overline{V_i},\overline{V_i} \cap B_\epsilon) = \sum_{W \subset \overline{V_i} \setminus V_i}  \frac{1}{b_{2e} \epsilon^{2e}} \frac{\eta (W,{\bf 1}_{\overline{V_i}})}{s_{2n-2e-1}} \int_{W \cap B_\epsilon} K_{2(d_W-e)} (x) dx + \hfill \cr
\hfill  \frac{1}{b_{2e} \epsilon^{2e}} \frac{1}{s_{2n-2e-1}} \int_{V_i \cap B_\epsilon} K_{2(d_i-e)} (x) dx, \qquad \cr
}$$
where $d_W$ denotes the dimension of the stratum $W$. Since 
$$\lim_{\epsilon \rightarrow 0} \frac{1}{b_{2e} \epsilon^{2e}} \Lambda_{2e} (\overline{V_i}, \overline{V_i} \cap B_\epsilon) =0,$$ and, by the induction hypothesis, 
$$\lim_{\epsilon \rightarrow 0} \frac{1}{b_{2e} \epsilon^{2e}}   \frac{1}{s_{2n-2e-1}} \int_{W \cap B_\epsilon} K_{2(d_W-e)} (x) dx=0,$$
we find that 
$$\lim_{\epsilon \rightarrow 0} \frac{1}{b_{2e}   \epsilon^{2e}}  \frac{1}{s_{2n-2e-1}} \int_{V_i \cap B_\epsilon} K_{2(d_i-e)}(x) dx=0.$$ 

The same argument works in order to prove that 
$$\lim_{\epsilon \rightarrow 0} \frac{1}{b_{2e} \epsilon^{2e}}  \frac{1}{s_{2n-2e-1}}\int_{V_i \cap B_\epsilon} K_{2(d_i-e)}(x) dx,$$ exists and is finite for $e \in \{d_0,\ldots,d_i \}$. \endproof

Let us focus now on the limits 
$$\lim_{\epsilon \rightarrow 0} \frac{1}{b_{2d_0} \epsilon^{2d_0}} \frac{1}{s_{2n-2d_0-1}} \int_{V_i \cap B_\epsilon} K_{2(d_i-d_0)} (x) dx,$$
for $i \in \{0,\ldots,q \}$. We already know that
$$\lim_{\epsilon \rightarrow 0} \frac{1}{b_{2d_0} \epsilon^{2d_0}} \frac{1}{s_{2n-2d_0-1}} \int_{V_0 \cap B_\epsilon} K_0(x) dx =1.$$

\begin{proposition}
For $i \in \{1,\ldots,q \}$, we have
$$\lim_{\epsilon \rightarrow 0} \frac{1}{b_{2d_0} \epsilon^{2d_0}} \frac{1}{s_{2n-2d_0-1}} \int_{V_i \cap B_\epsilon} K_{2(d_i-d_0)}(x) dx=0.$$
\end{proposition}
\proof By Theorem \ref{CurvAndPolarC}, we know that
$$\lim_{\epsilon \rightarrow 0} \frac{\Lambda_{2d_0}(X,X \cap B_\epsilon)}{b_{2d_0} \epsilon^{2d_0}} = 1 -\sigma_{2d_0 +1}(X,0).$$
Let $H \in G_{2n}^{2n-2d_0}$ be a generic linear space that intersects $V_0$ transversally at $0$. We can choose the scalar product in $\mathbb{C}^n= \mathbb{R}^{2n}$ in such a way that $H= T_0 V_0 ^\perp$. Since $0$ is not an exceptional point for $X$ is complex analytic, for almost all $v \in S_H^\perp$, the lower half-link 
$X \cap \{v^* = -\delta \} \cap B_\epsilon \cap H$, $0 < \delta \ll \epsilon \ll 1$, does not depend on the scalar product, on $H$ and on $v$ and we have
$$\chi \Big( X \cap \{v^* = -\delta \} \cap B_\epsilon \cap H \Big)= \chi(\mathcal{L}_{V_0}).$$
Since $H$ and $v$ are generic, we see that
$$ \lim_{\epsilon \rightarrow 0} \frac{\Lambda_{2d_0} (X, X \cap B_\epsilon)}{b_{2d_0} \epsilon^{2d_0}} = 1-\chi(\mathcal{L}_{V_0}) = \eta(V_0,{\bf 1}_X).$$
We can prove the proposition by induction on the depth of the stratum. Let $W$ be a stratum of depth equal to $1$. We have
$$\displaylines{
\quad \lim_{\epsilon \rightarrow 0} \frac{\Lambda_{2d_0} (\overline{W}, \overline{W} \cap B_\epsilon)}{b_{2d_0} \epsilon^{2d_0}} =  \hfill \cr
 \quad \quad \quad \eta(W,{\bf 1}_{\overline{W}}) \lim_{\epsilon \rightarrow 0} \frac{1}{b_{2d_0} \epsilon^{2d_0} } \frac{1}{s_{2n-2d_0-1}} \int_{W \cap B_\epsilon} K_{2(d_W-d_0)}(x) dx + \hfill \cr
 \hfill \eta(V_0,{\bf 1}_{\overline{W}}) \lim_{\epsilon \rightarrow 0} \frac{1}{b_{2d_0} \epsilon^{2d_0}} \frac{1}{s_{2n -2d_0-1}} \int_{V_0 \cap B_\epsilon} K_0(x) dx. \quad \cr
 }$$
By the above equality applied to $X=\overline{W}$ and the remark before the proposition, we find that
$$\eta(V_0,{\bf 1}_{\overline{W}} )= \lim_{\epsilon \rightarrow 0} \frac{1}{b_{2d_0} \epsilon^{2d_0}} \frac{1}{s_{2n-2d_0-1}} \int_{W \cap B_\epsilon} K_{2(d_W-d_0)}(x) dx + \eta(V_0,{\bf 1}_{\overline{W}} ).$$
Hence the result is true for any stratum of depth $1$. 

Let us fix now a stratum $V_i$. By induction, we can assume that the result is valid for any stratum $W$ included in $\overline{V_i} \setminus V_i$. We have
$$\displaylines{
\quad \lim_{\epsilon \rightarrow 0} \frac{\Lambda_{2d_0} (\overline{V_i}, \overline{V_i} \cap B_\epsilon)}{b_{2d_0} \epsilon^{2d_0}} = \eta(V_0,{\bf 1}_{\overline{V_i}} )= \hfill \cr 
\quad \quad \lim_{\epsilon \rightarrow 0} \frac{1}{b_{2d_0} \epsilon^{2d_0}} \frac{1}{s_{2n-2d_0-1}} \int_{V_i \cap B_\epsilon} K_{2(d_i-d_0)}(x) dx + \hfill \cr
\quad \quad  \quad \sum_{W \subset \overline{V_i} \setminus V_i \atop W \not= V_0} \eta (W,{\bf 1}_{\overline{V_i}}) \lim_{\epsilon \rightarrow 0}  \frac{1}{b_{2d_0} \epsilon^{2d_0}} \frac{1}{s_{2n-2d_0-1}} \int_{W \cap B_\epsilon} K_{2(d_W-d_0)}(x) dx + \hfill \cr
\hfill \eta(V_0,{\bf 1}_{\overline{V_i}} ) 
 \frac{1}{b_{2d_0} \epsilon^{2d_0}} \frac{1}{s_{2n-2d_0-1}} \int_{V_0 \cap B_\epsilon} K_0(x) dx . \quad \cr
 }$$
 Using the induction hypothesis, we get the result for the stratum $V_i$. \endproof
 
 Now we can state the main results of this section.
 \begin{theorem}\label{MainTheorem}
 Let $\phi : X \rightarrow \mathbb{Z}$ be a constructible function. We have
 $$\displaylines{
 \quad \phi (0) = \eta(V_0,\phi) + \hfill \cr
\quad \quad \sum_{i=1}^q \left[ \sum_{e=d_0+1}^{d_i} \lim_{\epsilon \rightarrow 0} \frac{1}{b_{2e} \epsilon^{2e}} \frac{1}{s_{2n-2e-1}} \int_{V_i \cap B_\epsilon} K_{2d_i -2e}(x) dx \right] \eta(V_i, \phi) = \hfill\cr
 \hfill  \eta(V_0,\phi) + 
 \sum_{i=1}^q \left[ \sum_{e=d_0+1}^{d_i} \lim_{\epsilon \rightarrow 0} \frac{1}{e ! b_{2e} \epsilon^{2e}} \int_{V_i \cap B_\epsilon}   {\rm ch}_{d_i-e}(V_i) \wedge \kappa(V_i)^e \right] \eta(V_i, \phi). \quad \cr 
 }$$
 \end{theorem}
 \proof Let us prove the first equality. For $\phi={\bf 1}_X$, we just  reformulate Corollary \ref{CorCurvAndLink2} using the previous two propositions. Let $V$ be a stratum of $X$. Then the formula is also true for $\phi={\bf 1}_{\overline{V}}$ because $\eta (V_i,{\bf 1}_{\overline{V}})=0$ if $V_i \nsubseteq \overline{V}$. Since both sides of the equality are linear in $\phi$, this gives the result. 
 
 The second equality is an application of Lemma \ref{LipschitzChern}. \endproof
 
 \begin{corollary}\label{Euler}
 Assume that $X$ is equidimensional. We have
 $$\displaylines{
 \qquad {\rm Eu}_X(0)= \sum_{e=d_0+1}^d \lim_{\epsilon \rightarrow 0} \frac{1}{b_{2e} \epsilon^{2e}} \frac{1}{s_{2n-2e-1}} \int_{X_{{\rm reg}} \cap B_\epsilon} K_{2(d-e)}(x) dx= \hfill \cr
 \hfill \sum_{e=d_0+1}^d \lim_{\epsilon \rightarrow 0} \frac{1}{e! b_{2e} \epsilon^{2e}}  \int_{X_{{\rm reg}} \cap B_\epsilon}  {\rm ch}_{d-e}(X_{{\rm reg}}) \wedge \kappa(X_{{\rm reg}})^e . \qquad \cr
 }$$
 \end{corollary}
 \proof In this situation $X_{{\rm reg}}$ is exactly the union of the top strata. We apply the previous theorem to $\phi={\rm Eu}_X$ and use the fact that $\eta(V_i, {\rm Eu}_X)=0$ if $V_i$ is not a top stratum and that $\eta(V_i, {\rm Eu}_X)=1$ if $V_i$ is a top stratum. \endproof
 
 \begin{remark}
 As mentioned in the introduction, we think that the above corollary can be proved by complex methods using Loeser's results \cite{Loeser}.  Moreover, motivated by  the relation between the Euler obstruction and the polar multiplicities established by L\^e and Teissier \cite{LeTeissierAnnals} (see also \cite{Massey}), we believe that the following 
 equality shoud be true:
 $$ \lim_{\epsilon \rightarrow 0} \frac{1}{e! b_{2e} \epsilon^{2e}}  \int_{X_{{\rm reg}} \cap B_\epsilon}  {\rm ch}_{d-e}(X_{{\rm reg}}) \wedge \kappa(X_{{\rm reg}})^e = (-1)^e m_0(\Gamma_e),$$
 where $m_0(\Gamma_e)$ is the multiplicity of the  polar variety of codimension $e$.
 \end{remark}
 As a corollary, we recover the local index formula of Brylinski, Dubson and Kashiwara \cite{BDK}, which in this general form is due to Sch\"urmann \cite{Schuermann} (Equality 5.40, p.294). Note that since here we use the real orientation in the definition of the Euler obstruction instead of the complex orientation,
 our formula differs from Sch\"urmann's equality, where each $\eta(V_i,\phi)$ is replaced by $(-1)^{d_i} \eta(V_i,\phi)$ (see remark 5.0.3, p.293 in \cite{Schuermann}).
 \begin{corollary}\label{IndiceBDK}
 Assume that $X$ is equidimensional. Let $\phi : X \rightarrow \mathbb{Z}$ be a constructible function. We have
 $$\phi(0)= \sum_{i=0}^q {\rm Eu}_{\overline{V_i}}(0) \eta(V_i,\phi).$$
 \end{corollary}
 \proof Apply the previous corollary to $X=\overline{V_i}$ and use the fact that ${\rm Eu}_{\overline{V_0}}(0)=1$. \endproof
 
 \section{Global Euler obstruction and curvatures}
 In this section, we give global versions of the results that we established in the previous section.
 
 We recall first the Gauss-Bonnet formula for closed semi-algebraic sets that we proved in \cite{DutertreGeoDedicata}. Let $X \subset \mathbb{R}^n$ be a closed semi-algebraic set of dimension $d$, equipped with a semi-algebraic Whitney stratification $\{ V_a \}_{a \in A}$. The Lipschitz-Killing measures 
 $$\Lambda_0(X,-),\ldots,\Lambda_n(X,-),$$
 are defined as in the  subanalytic case.  
 
 In  \cite{DutertreGeoDedicata}, Theorem 3.3, we proved the following Gauss-Bonnet formula.
 \begin{theorem}\label{GBSemiAlg}
Let $X\subset \mathbb{R}^n$ be a closed semi-algebraic set of dimension $d$. We have
$$\chi(X)=\sum_{k=0}^d    \lim_{R \rightarrow +\infty} \frac{\Lambda_k(X,X \cap B_R)}{b_k R^k}.$$
\end{theorem}
$\hfill \Box$

Next we apply this formula when $X \subset \mathbb{C}^n$ is a complex algebraic set. We write $X=\cup_{i=0}^q V_i$. We assume that the stratum $V_0,\ldots,V_q$ are connected and that the analytic sets $\overline{V_0},\ldots,\overline{V_q}$ are reduced. We observe that $V_0=\overline{V_0}$. We set $d_i= {\rm dim} V_i$ for $i \in \{0,\ldots,q\}$. %We call a {\bf top stratum} a stratum not contained in the frontier of any other stratum. 

Let $\phi : X \rightarrow \mathbb{Z}$ be a constructible function with respect to this stratification. In this situation and using the same arguments as in the previous section, the above Gauss-Bonnet formula leads to the following theorem.
 \begin{theorem}\label{GlobalMainTheorem}
 Let $\phi : X \rightarrow \mathbb{Z}$ be a constructible function. We have
 $$\displaylines{
 \quad \chi(X,\phi) = \hfill \cr
\quad \quad \sum_{i=0}^q \left[ \sum_{e=0}^{d_i} \lim_{R \rightarrow +\infty} \frac{1}{b_{2e} R^{2e}} \frac{1}{s_{2n-2e-1}} \int_{V_i \cap B_R} K_{2d_i -2e}(x) dx \right] \eta(V_i, \phi) = \hfill\cr
\hfill  \sum_{i=0}^q \left[ \sum_{e=0}^{d_i} \lim_{R \rightarrow +\infty} \frac{1}{e ! b_{2e} R^{2e}} \int_{V_i \cap B_R}   {\rm ch}_{d_i-e}(V_i) \wedge \kappa(V_i)^e \right] \eta(V_i, \phi). \quad \cr 
 }$$
 \end{theorem}
$\hfill \Box$

Now we assume that $X$ is equidimensional. In this case, Seade, Tib\u{a}r and Verjovsky introduced a global analogous of the Euler obstruction called the global Euler obstruction and denoted by ${\rm Eu}(X)$ (see \cite{SeadeTibarVerjovsky}, Definition 3.3). They also gave a global version of the L\^e-Teissier formula relating ${\rm Eu}(X)$ to global polar invariants (\cite{SeadeTibarVerjovsky}, Theorem 3.4).  This result was later generalized by Sch\"urmann and Tib\u{a}r in \cite {SchuermannTibar}  using the language of constructible functions and Mac-Pherson cycles. 

When we apply the above theorem to $\phi={\rm Eu}_X$, we get the following Gauss-Bonnet formula for the global Euler obstruction.
\begin{corollary}\label{GlobalEuler}
 Assume that $X$ is equidimensional. We have
 $$\displaylines{
 \qquad {\rm Eu}(X)= \sum_{e=0}^d \lim_{R \rightarrow +\infty} \frac{1}{b_{2e} R^{2e}} \frac{1}{s_{2n-2e-1}} \int_{X_{{\rm reg}} \cap B_R} K_{2(d-e)}(x) dx= \hfill \cr
 \hfill \sum_{e=0}^d \lim_{R \rightarrow +\infty} \frac{1}{e! b_{2e}R^{2e}}  \int_{X_{{\rm reg}} \cap B_R}  {\rm ch}_{d-e}(X_{{\rm reg}}) \wedge \kappa(X_{{\rm reg}})^e . \qquad \cr
 }$$
 \end{corollary}
$\hfill \Box$

We remark that as in the local case, each member of this last sum should be equal up to sign to a global polar invariant. As a corollary, we obtain a global version of the Brylinski-Dubson-Kashiwara formula.
\begin{corollary}\label{GlobalIndiceBDK}
 Assume that $X$ is equidimensional. Let $\phi : X \rightarrow \mathbb{Z}$ be a constructible function. We have
 $$\chi(X,\phi)= \sum_{i=0}^q {\rm Eu}(\overline{V_i}) \eta(V_i,\phi).$$
 \end{corollary}
$\hfill \Box$

We note that for $\phi={\bf 1}_X$, this formula was proved by Tib\u{a}r \cite{TibarEuler}.

\section{Euler obstruction and the Gauss-Bonnet measure}
In this section, using our results on the Euler obstruction, we give a positive answer to a question of Fu \cite{FuJDiffGeo} on the Euler obstruction and the Gauss-Bonnet measure. 
 
 Let us go back before to the subanalytic case. Let $X \subset \mathbb{R}^n$ be a closed subanalytic set such that $0 \in X$. Since we work in a neighborhood of $0$, we can assume that $X$ is equipped with a finite Whitney subanalytic stratification $X=\cup_{i=0}^q V_i$, where $0$ belongs to the closure of each stratum $V_i$ and $V_0$ is the stratum containing $0$. For $\epsilon >0$ sufficiently small, $B_\epsilon$ intersects $X$ transversally and so $X \cap B_\epsilon$ admits the following Whitney stratification:
 $$X \cap B_\epsilon = \cup_{i=0}^q V_i \cap \mathring{B_\epsilon} \bigcup \cup_{i=0}^q V_i \cap S_\epsilon,$$
 (note that $V_0 \cap S_\epsilon = \emptyset$ if $V_0=\{0\})$. By the Gauss-Bonnet theorem of Fu \cite{Fu94} and Broecker and Kuppe \cite{BroeckerKuppe}, we have
 $$\displaylines{
 \qquad \chi(X \cap B_\epsilon)= \Lambda_0 (X \cap B_\epsilon, X \cap B_\epsilon) = \Lambda_0(X \cap B_\epsilon, X \cap \mathring{B_\epsilon}) + \hfill \cr
 \hfill  \Lambda_0 (X \cap B_\epsilon, X \cap S_\epsilon) =
 \Lambda_0(X , X \cap B_\epsilon) + \Lambda_0(X \cap B_\epsilon, X \cap S_\epsilon). \qquad \cr
 }$$
 Since $\lim_{\epsilon \rightarrow 0} \chi(X \cap B_\epsilon) = 1$ and, by Theorem 5.1 in \cite{DutertreGeoDedicata}, $\lim_{\epsilon \rightarrow 0} \Lambda_0(X,X \cap B_\epsilon)$ exists and is finite, we find that $\lim_{\epsilon \rightarrow 0} \Lambda_0(X \cap B_\epsilon, X \cap S_\epsilon)$ exists and is finite. In the sequel, we will give a characterization of this limit in terms of indices of critical points on $X \cap S_\epsilon$ of generic linear functions. We will apply Theorem 3.1 and Lemma 2.1 of \cite{DutertreManuscripta}. 
 Let us recall first the definition of the index of an isolated stratified critical point.
\begin{definition}
Let $Z \subset \mathbb{R}^n$ be a closed subanalytic set, equipped with a Whitney stratification. Let $p \in Z$ be an isolated critical point of a subanalytic function $g : Z \rightarrow \mathbb{R}$, which is the restriction to $Z$ of a $C^2$-subanalytic function $G$. We define the index of $g$ at $p$ as follows:
$${\rm ind}(g,Z,p)= 1 - \chi \big( Z \cap \{ g = g(p)-\delta \} \cap B_\epsilon(p) \big),$$
where $0< \delta \ll \epsilon \ll 1$.
\end{definition}
The following lemma  is necessary in order to apply Lemma 2.1 of \cite{DutertreManuscripta}.
 \begin{lemma}\label{FirstLemmaSection}
 Let $V$ be a stratum of $X$ such that $0 < {\rm dim} V < n$. For almost all $v \in S^{n-1}$, there exists $\epsilon_v >0$ such that for $0 < \epsilon \le \epsilon_v$, $v^*_{\vert S_\epsilon}$ is a submersion at any point of $S_\epsilon \cap V$, and so at any critical point of $v^*_{\vert S_\epsilon \cap V}$.
 \end{lemma}
 \proof By Lemma 5.4 in \cite{DutertreProcTrotman} applied to $H= \mathbb{R}^n$, we know that for almost all $v$ in $S^{n-1}$, the line $l_v$ generated by $v$ intersects $V$ transversally. This implies that the intersection $l_v \cap V$ has dimension at most $0$ and so $l_v$ does not intersect $V$ in a small punctured neighborhood of the origin. \endproof
 
\begin{lemma}\label{SecondLemmaSection}
For almost all $v$ in $S^{n-1}$, $v^*_{\vert X}$ has an isolated critical point at the origin.
\end{lemma}
\proof See \cite{DutertreProcTrotman}, Corollary 4.2. \endproof

Let $(\epsilon_s)_{s \in \mathbb{N}}$ be a sequence of positive real number such that $\lim_{s \rightarrow + \infty} \epsilon_s =0$. There exists $s_0$ such that for $s \ge s_0$, $B_{\epsilon_s}$ intersects $X$ transversally and $B_{\epsilon_s} \cap X$ is naturally Whitney stratified. Let $s$ be such that $s \ge s_0$. By \cite{BroeckerKuppe} Lemma 3.5, for almost all $v$ in $S^{n-1}$, $v^*_{\vert X \cap B_{\epsilon_s}}$ is a Morse function. Since a countable union of sets of measure zero has measure zero, for almost all $v$ in $S^{n-1}$ the function $v^*_{\vert X \cap B_{\epsilon_s}}$ is a Morse function for any $s \ge s_0$.  

Let us fix a generic $v$ in $S^{n-1}$ which satisfies this condition and the conditions of Lemma \ref{FirstLemmaSection} and \ref{SecondLemmaSection}. There exists $s_v$ such that for any $s \ge s_v$, $v^*_{\vert X \cap B_{\epsilon_s}}$ is a Morse function, $0$ is its only critical point lying in $\mathring{B_{\epsilon_s}}$ and $v^*_{\vert S_{\epsilon_s}}$ is a submersion at any critical point of $v^*_{\vert X \cap S_{\epsilon_s}}$. Let $p$ be a critical point of $v^*_{\vert X \cap B_{\epsilon_s}}$ in $S_{\epsilon_s}$ and let $V$ be the stratum that contains it. Since $v^*_{\vert X}$ has no critical point on $V \setminus \{0\}$, there exists $\lambda(p) \not= 0$ such that 
$$ \nabla (v^*_{\vert S})(p)=\lambda (p) \nabla(\omega_{\vert S})(p),$$
where $\omega$ is the euclidian distance function. We say that $p$ is outwards-pointing (resp. inwards-pointing) if $\lambda (p) >0$ (resp. $\lambda (p) <0$).  

We can apply Lemma 2.1 of \cite{DutertreManuscripta} and find
that ${\rm ind}(v^*, X \cap B_{\epsilon_s}, p)=0$ if $p$ is outwards-pointing and 
$${\rm ind}(v^*, X \cap B_{\epsilon_s}, p)={\rm ind}(v^*, X \cap S_{\epsilon_s},p),$$
if $p$ is inwards-pointing. 

Since $v^*_{\vert X \cap B_{\epsilon_s}}$ is a Morse function, $v^*_{\vert X \cap S_{\epsilon_s}}$ is a Morse function as well. Therefore we can write 
$${\rm ind}(v^*,X \cap S_{\epsilon_s},p)= (-1)^{\sigma(p)} \cdot {\rm ind}_{nor} (v^*,X \cap S_{\epsilon_s},p),$$
where $\sigma(p)$ is the Morse index of $v^*_{\vert V \cap S_{\epsilon_s}}$ at $p$ and where ${\rm ind}_{nor} (v, X \cap S_{\epsilon_s},p)$ is the normal Morse index at $p$. It is defined as follows:
$${\rm ind}_{nor} (v, X \cap S_{\epsilon_s},p)= 1-\chi \Big( X \cap S_{\epsilon_s} \cap N_p \cap B_{\nu}(p) \cap \{ v^* = v^* (p)-\delta \} \Big),$$
where $0 < \delta \ll \nu \ll 1$ and $N_p$ is a normal slice to the stratum $V \cap S_{\epsilon_s}$ such that $N_p \cap V \cap S_{\epsilon_s} = \{ p\}$. % and $B_{\nu}(p)$ is the closed ball of radius $\nu$ centered at $p$. 
Let us denote by $\mathcal{I}_{v,s}$ the set of inwards-pointing critical points of $v^*_{\vert X \cap B_{\epsilon_s}}$. By Theorem 3.1 in \cite{DutertreManuscripta}, we have that for $s \ge s_0$
$$\chi(X \cap B_{\epsilon_s})=  {\rm ind}(v^*,X,0) + \sum_{p \in \mathcal{I}_{v,s} } (-1)^{\sigma(p)} \cdot {\rm ind}_{nor} (v^*, X \cap S_{\epsilon_s}, p). \eqno(*)$$
Since $\lim_{s \rightarrow +\infty} \chi(X \cap B_{\epsilon_s})=1$, we see that $$\lim_{s \rightarrow + \infty} \sum_{p \in \mathcal{I}_{v,s}} (-1)^{\sigma(p)} \cdot {\rm ind}_{nor}(v^*, X \cap S_{\epsilon_s},p),$$ exists and is equal to $1-{\rm ind}(v^*,X,0)$. Taking the mean-value on $S^{n-1}$ in this last equality and passing to the limit, we get that
$$\displaylines{
\qquad 1-\lim_{s \rightarrow + \infty} \Lambda_0(X,X \cap B_{\epsilon_s}) = \hfill \cr
\hfill  \frac{1}{s_{n-1}} \int_{S^{n-1}} \lim_{s \rightarrow + \infty} \sum_{p \in \mathcal{I}_{v,s}} (-1)^{\sigma(p)} \cdot {\rm ind}(v^*, X \cap S_{\epsilon_s}, p) dv. \qquad \cr
}$$
\begin{proposition}\label{ThirdPropSection}
We have
$$\displaylines{
\qquad \lim_{\epsilon \rightarrow 0} \Lambda_0 ( X \cap B_\epsilon, X \cap S_\epsilon) = \lim_{s \rightarrow + \infty} \Lambda_0(X \cap B_{\epsilon_s},  X \cap S_{\epsilon_s}) = \hfill \cr
\hfill  \frac{1}{s_{n-1}} \int_{S^{n-1}} \lim_{s \rightarrow + \infty} \sum_{p \in \mathcal{I}_{v,s}} (-1)^{\sigma(p)} \cdot {\rm ind}_{nor}(v^*, X \cap S_{\epsilon_s}, p) dv, \qquad \cr
}$$
for any sequence $(\epsilon_s)_{s \in \mathbb{N}}$ of positive real numbers such that $\lim_{s \rightarrow + \infty} \epsilon_s =0$.
\end{proposition}
$\hfill \Box$

We believe that this equality has its own interest in the subanalytic case because it gives a topological description of $\lim_{\epsilon \rightarrow 0} \Lambda_0(X \cap B_\epsilon, X \cap S_\epsilon)$. In the sequel, we will refine it when $X$ is a complex analytic set.

Let $(X,0) \subset (\mathbb{C}^n,0)$ be the germ of an analytic complex variety. We keep the notations of Section 4. We consider a sequence $(\epsilon_s)_{s \in \mathbb{N}}$ of positive real numbers tending to $0$ and a vector $v$ in $S^{2n-1}$ generic as above and we choose $s \ge s_v$. Let $p \in V \cap S_{\epsilon_s}$ be an inwards-pointing critical point of $v^*_{\vert X \cap B_{\epsilon_s}}$. Let $N_p$ be a normal slice to $V \cap S_{\epsilon_s}$ such that $N_p \cap (V \cap S_{\epsilon_s}) =\{p\}$. Then $N_p \cap S_{\epsilon_s}$ is a normal slice to $V$ such that $(N_p \cap S_{\epsilon_s}) \cap V =\{p\}$. Moreover, since the form $v^*_{\vert N_p}$ is non-degenerate for $X \cap S_{\epsilon_s} \cap N_p$ at $p$, the form $v^*_{\vert N_p \cap S_{\epsilon_s}}$ is non-degenerate for $X \cap S_{\epsilon_s} \cap N_p$ at $p$ as well. Therefore the set
$$X \cap N_p \cap S_{\epsilon_s} \cap B_{\nu}(p) \cap \{v^*=v^*(p)-\delta \},$$
is the lower half-link of $v^*_{\vert X}$ and, as already explained in Section 3, its Euler characteristic is equal to $\chi(\mathcal{L}_V)$. Equality $(*)$ becomes
$$\chi(X \cap B_{\epsilon_s}) ={\rm ind}(v^*,X,0) + \sum_{i=0}^q \eta(V_i,{\bf 1}_X) \sum_{p \in \mathcal{I}^i_{v,s}} (-1)^{\sigma(p)},$$ 
where $\mathcal{I}^i_{v,s}$ is the set of inwards-pointing critical points of $v^*_{\vert X \cap B_{\epsilon_s}}$ in $V_i \cap S_{\epsilon_s}$. Applied to $X=\overline{V_0}$, this gives that $\lim_{s \rightarrow + \infty} \sum_{p \in \mathcal{I}^0_{v,s}} (-1)^{\sigma(p)} $ exists, is finite and does not depend on the choice of the sequence. Note that if ${\rm dim} V_0  =0$ then $\mathcal{I}_{v,s}^0$ is empty and $\lim_{s \rightarrow + \infty} \sum_{p \in \mathcal{I}^0_{v,s}} (-1)^{\sigma(p)}=0 $. 
Applied to $X=\overline{V_j}$, where $V_j$ is a stratum of depth $1$, it gives that $\lim_{s \rightarrow + \infty} \sum_{p \in \mathcal{I}^j_{v,s}} (-1)^{\sigma(p)}$ exists, is finite and does not depend on the choice of the sequence. By induction on the depth of the stratum, we see that for $i \in \{0,\ldots,q\}$, $\lim_{s \rightarrow + \infty} \sum_{p \in \mathcal{I}^i_{v,s}} (-1)^{\sigma(p)}$ exists, is finite and does not depend on the choice of the sequence.
Proposition \ref{ThirdPropSection} becomes
\begin{proposition}
We have 
$$\lim_{\epsilon \rightarrow 0} \Lambda_0 ( X \cap B_\epsilon, X \cap S_\epsilon) = \lim_{s \rightarrow + \infty} \Lambda_0(X \cap B_{\epsilon_s},  X \cap S_{\epsilon_s}) =$$
$$ \frac{1}{s_{2n-1}} \int_{S^{2n-1}} \sum_{i=0}^q \eta(V_i,{\bf 1}_X) \cdot \lim_{s \rightarrow + \infty} \sum_{p \in \mathcal{I}^i_{v,s} }(-1)^{\sigma(p) }dv =$$
$$\sum_{i=0}^q \eta(V_i,{\bf 1}_X)  \frac{1}{s_{2n-1}} \int_{S^{2n-1}}  \lim_{s \rightarrow + \infty} \sum_{p \in \mathcal{I}^i_{v,s}} (-1)^{\sigma(p) }dv ,$$
for any sequence $(\epsilon_s)_{s \in \mathbb{N}}$ of positive real numbers such that $\lim_{s \rightarrow + \infty} \epsilon_s =0$.
\end{proposition}
$\hfill \Box$

\begin{corollary}
For $i \in \{0,\ldots,q\}$, $\lim_{\epsilon \rightarrow 0} \Lambda_0(\overline{V_i} \cap B_\epsilon, V_i \cap S_\epsilon)$ exists and is finite. Furthermore, we have
$$\lim_{\epsilon \rightarrow 0} \Lambda_0(X\cap B_\epsilon, X \cap S_\epsilon) = \sum_{i=0}^q \eta(V_i,{\bf 1}_X) \cdot \lim_{\epsilon \rightarrow 0} \Lambda_0(\overline{V_i} \cap B_\epsilon, V_i \cap S_\epsilon).$$
\end{corollary}
\proof We remark first that $\lim_{\epsilon \rightarrow 0} \Lambda_0 (\overline{V_0} \cap B_\epsilon, V_0 \cap S_\epsilon)=0$ if ${\rm dim} V_0=0$. 

Let $(\epsilon_s)_{s \in \mathbb{N}}$ be a sequence of positive real numbers such that $\lim_{s \rightarrow + \infty} \epsilon_s =0$. We have
$$\Lambda_0(\overline{V_i} \cap B_{\epsilon_s}, V_i \cap S_{\epsilon_s})=\frac{1}{s_{2n-1}} \int_{S^{2n-1}} \sum_{q \in V_i \cap S_{\epsilon_s}} {\rm ind}(v^*, \overline{V_i} \cap B_{\epsilon_s},q) dv.$$
Taking the limit as $s$ tends to $+ \infty$, we find
$$\lim_{s \rightarrow + \infty} \Lambda_0(\overline{V_i} \cap B_{\epsilon_s}, V_i \cap S_{\epsilon_s}) = \frac{1}{s_{2n-1}} \int_{S^{2n-1}} \lim_{s \rightarrow + \infty} \sum_{q \in V_i \cap S_{\epsilon_s}} {\rm ind} (v^*, \overline{V_i} \cap B_{\epsilon_s}, q) dv.$$
By our previous study, we know that for $v$ generic,
$$\lim_{s \rightarrow + \infty} \sum_{q \in V_i \cap S_{\epsilon_s}} {\rm ind}(v^*, \overline{V_i} \cap B_{\epsilon_s},q)= \lim_{s \rightarrow +\infty} \sum_{p \in \mathcal{I}^i_{v,s}} (-1)^{\sigma(p)},$$
and so, this limit exists and is finite. Therefore, we find that
$$\lim_{s \rightarrow + \infty} \Lambda_0(\overline{V_i} \cap B_{\epsilon_s}, V_i \cap S_{\epsilon_s}),$$ exists and is equal to
$$\frac{1}{s_{2n-1}} \int_{S^{2n-1}} \lim_{s \rightarrow + \infty} \sum_{p \in \mathcal{I}^i_{v,s}} (-1)^{\sigma(p)} dv.$$
Since this last integral does not depend on the choice of the sequence, we can conclude that $\lim_{\epsilon \rightarrow 0} \Lambda_0(\overline{V_i} \cap B_\epsilon, V_i \cap S_\epsilon)$ exists and is equal to
$$\frac{1}{s_{2n-1}} \int_{S^{2n-1}} \lim_{s \rightarrow + \infty} \sum_{p \in \mathcal{I}^i_{v,s} }(-1)^{\sigma(p)} dv.$$
The previous proposition enables us to finish the proof. \endproof

\begin{theorem}
If ${\rm dim} V_0 >0$ then $\lim_{\epsilon \rightarrow 0} \Lambda_0(\overline{V_0} \cap B_\epsilon, V_0 \cap S_\epsilon)=1$. Furthermore,  for $i \in \{1,\ldots,q\}$, we have
$$\displaylines{
\quad \lim_{\epsilon \rightarrow 0} \Lambda_0(\overline{V_i} \cap B_\epsilon, V_i \cap S_\epsilon) = 
\sum_{e=d_0+1}^{d_i} \lim_{\epsilon \rightarrow 0} \frac{1}{b_{2e} \epsilon^{2e}} \frac{1}{s_{2n-2e-1}} \int_{V_i \cap B_\epsilon} K_{2d_i-2e}(x) dx= \hfill \cr
\hfill \sum_{e=d_0+1}^{d_i}  \lim_{\epsilon \rightarrow 0} \frac{1}{e ! b_{2e} \epsilon^{2e}} \int_{V_i \cap B_\epsilon}   {\rm ch}_{d_i-e}(V_i) \wedge \kappa(V_i)^e . \quad \cr
}$$
\end{theorem}
\proof To prove the first equality, we use the fact that 
$$\lim_{\epsilon \rightarrow 0} \Lambda_0(\overline{V_0} \cap B_\epsilon, V_0 \cap S_\epsilon)=1-\lim_{\epsilon \rightarrow 0} \Lambda_0(\overline{V_0},V_0 \cap B_\epsilon),$$
and the fact that $\lim_{\epsilon \rightarrow 0} \Lambda_0 (\overline{V_0}, V_0 \cap B_\epsilon)=0$ since $V_0$ is smooth. 

We prove the second equality by induction on the depth of the stratum. We assume first that ${\rm dim} V_0 >0$. Let $W$ be a stratum of depth $1$. By the previous corollary, we have
$$\lim_{\epsilon \rightarrow 0} \Lambda_0(\overline{W} \cap B_\epsilon, \overline{W} \cap S_\epsilon) = \eta (V_0,{\bf 1}_{\overline{W}}) + \lim_{\epsilon \rightarrow 0} \Lambda_0(\overline{W} \cap B_\epsilon, W \cap S_\epsilon).$$
But we also know that
$$\lim_{\epsilon \rightarrow 0} \Lambda_0(\overline{W} \cap B_\epsilon, \overline{W} \cap S_\epsilon) = 1-\lim_{\epsilon \rightarrow 0} \Lambda_0(\overline{W}  , \overline{W} \cap B_\epsilon).$$
By Corollaries \ref{CorCurvAndLink} and  \ref{CorCurvAndLink2} and the description of the limits $\lim_{\epsilon \rightarrow 0} \frac{\Lambda_k (\overline{W}, \overline{W} \cap B_\epsilon)}{b_k \epsilon ^k}$ given in Section 4, we get
$$\displaylines{
\qquad 1-\lim_{\epsilon \rightarrow 0} \Lambda_0 (\overline{W}, \overline{W} \cap B_\epsilon) = \eta(V_0,{\bf 1}_{\overline{W}} )+ \hfill \cr
\hfill \sum_{e=d_0+1}^{d_W} \lim_{\epsilon \rightarrow 0} \frac{1}{b_{2e} \epsilon^{2e}} \frac{1}{s_{2n-2e-1}} \int_{W \cap B_\epsilon} K_{2d_W-2e}(x) dx, \qquad \cr
}$$
where $d_W$ is the dimension of $W$. Comparing these relations, we obtain the result for a stratum of depth $1$. Let us prove the result for a stratum $V_i$, $i \in \{1,\ldots,q \}$. We have
$$\displaylines{
\qquad 1-\lim_{\epsilon \rightarrow 0} \Lambda_0(\overline{V_i}, \overline{V_i} \cap B_\epsilon) = \lim_{\epsilon \rightarrow 0} \Lambda_0(\overline{V_i} \cap B_\epsilon, \overline{V_i} \cap S_\epsilon)= \hfill \cr
\hfill \sum_{W \subset \overline{V_i} \setminus V_i  } \eta(W,{\bf 1}_{\overline{V_i}} ) \lim_{\epsilon \rightarrow 0} \Lambda_0 (\overline{W} \cap B_\epsilon, W \cap S_\epsilon) +  \lim_{\epsilon \rightarrow 0} \Lambda_0 (\overline{V_i} \cap B_\epsilon, V_i \cap S_\epsilon)  , \qquad \cr
}$$
and 
$$\displaylines{
\quad 1-\lim_{\epsilon \rightarrow 0} \Lambda_0(\overline{V_i}, \overline{V_i} \cap B_\epsilon) =  \lim_{\epsilon \rightarrow 0} \frac{1}{b_{2e} \epsilon^{2e}} \frac{1}{s_{2n-2e-1}} \int_{V_i \cap B_\epsilon} K_{2d_i-2e}(x)d x + \hfill \cr
\quad \quad  \sum_{W \subset \overline{V_i} \setminus V_i \atop W \not= V_0} \eta(W,{\bf 1}_{\overline{V_i}} ) \cdot
\sum_{e=d_0+1}^{d_W} \lim_{\epsilon \rightarrow 0} \frac{1}{b_{2e} \epsilon^{2e}} \frac{1}{s_{2n-2e-1}} \int_{W \cap B_\epsilon} K_{2d_W-2e}(x)d x  +\hfill \cr
\hfill  \eta(V_0,\bf{1}_{\overline{V_i}} ). \qquad \cr
}$$
The result is obtained applying the induction hypothesis.

If ${\rm dim} V_0=0$ and $W$ is a stratum of depth $1$, then we have
$$\lim_{\epsilon \rightarrow 0} \Lambda_0(\overline{W} \cap B_\epsilon, \overline{W} \cap S_\epsilon) =  \lim_{\epsilon \rightarrow 0} \Lambda_0(\overline{W} \cap B_\epsilon, W \cap S_\epsilon),$$
and
$$\lim_{\epsilon \rightarrow 0} \Lambda_0(\overline{W} \cap B_\epsilon, \overline{W} \cap S_\epsilon) = 1-\lim_{\epsilon \rightarrow 0} \Lambda_0(\overline{W}  , \overline{W} \cap B_\epsilon).$$
But, in this situation, Corollary \ref{CorCurvAndLink2} becomes
$$
1-\lim_{\epsilon \rightarrow 0} \Lambda_0 (\overline{W}, \overline{W} \cap B_\epsilon) =  
 \sum_{e= 1}^{d_W} \lim_{\epsilon \rightarrow 0} \frac{1}{b_{2e} \epsilon^{2e}} \frac{1}{s_{2n-2e-1}} \int_{W \cap B_\epsilon} K_{2d_W-2e}(x) dx, $$
because the stratum $\{0\}$ has no contribution in the computation of the curvatures $\Lambda_k (\overline{W}, \overline{W} \cap B_\epsilon)$, $k \ge 1$. This gives the result for a stratum of depth $1$. Let $V_i$ be a stratum, $i \in \{1,\ldots,q\}$. We have
$$\displaylines{
\qquad 1-\lim_{\epsilon \rightarrow 0} \Lambda_0(\overline{V_i}, \overline{V_i} \cap B_\epsilon) = \lim_{\epsilon \rightarrow 0} \Lambda_0(\overline{V_i} \cap B_\epsilon, \overline{V_i} \cap S_\epsilon)= \hfill \cr
\hfill \sum_{W \subset \overline{V_i} \setminus V_i \atop W \not= V_0 } \eta(W,{\bf 1}_{\overline{V_i}} ) \lim_{\epsilon \rightarrow 0} \Lambda_0 (\overline{W} \cap B_\epsilon, W \cap S_\epsilon) +  \lim_{\epsilon \rightarrow 0} \Lambda_0 (\overline{V_i} \cap B_\epsilon, V_i \cap S_\epsilon)  , \qquad \cr
}$$
and, because the stratum $V_0$ is zero-dimensional,
$$\displaylines{
\quad 1-\lim_{\epsilon \rightarrow 0} \Lambda_0(\overline{V_i}, \overline{V_i} \cap B_\epsilon) =  \lim_{\epsilon \rightarrow 0} \frac{1}{b_{2e} \epsilon^{2e}} \frac{1}{s_{2n-2e-1}} \int_{V_i \cap B_\epsilon} K_{2d_i-2e}(x)d x + \hfill \cr
\hfill  \quad  \sum_{W \subset \overline{V_i} \setminus V_i \atop W \not= V_0} \eta(W,{\bf 1}_{\overline{V_i}} ) \cdot
\sum_{e=1}^{d_W} \lim_{\epsilon \rightarrow 0} \frac{1}{b_{2e} \epsilon^{2e}} \frac{1}{s_{2n-2e-1}} \int_{W \cap B_\epsilon} K_{2d_W-2e}(x)d x . \quad \cr
}$$
We obtain the result comparing these two equalities and applying the induction hypothesis.
\endproof

The following corollary gives a positive answer to Fu's question.
\begin{corollary}\label{FuQuestion}
Assume that $X$ is equidimensional. We have
$${\rm Eu}_X(0) = \lim_{\epsilon \rightarrow 0} \Lambda_0 (X \cap B_\epsilon, X_{{\rm reg}} \cap S_\epsilon).$$
\end{corollary}
\proof We apply Corollary \ref{Euler} and the fact that $X_{{\rm reg}}$ is exactly the union of the top strata. \endproof

\end{document}